\documentclass[a4paper,fleqn,leqno]{article}
\def\ladate{26 mai 2011. Nouvelle introduction: 6 janvier 2012.}

\usepackage[latin1]{inputenc}

\usepackage{geometry}
\usepackage[T1]{fontenc}

\usepackage[german,french]{babel}
\usepackage{amsmath,amsthm,amssymb}
\allowdisplaybreaks

\usepackage{enumerate} 

\usepackage{mathtools}
\DeclarePairedDelimiter\OO{\rbrack}{\lbrack}
\DeclarePairedDelimiter\OF{\rbrack}{\rbrack}
\DeclarePairedDelimiter\FO{\lbrack}{\lbrack}

\usepackage[pdfencoding=pdfdoc]{hyperref}
\hypersetup{%
linkcolor={blue},
citecolor={red},
colorlinks,%
bookmarksopen,%
pdfdisplaydoctitle,%
pdfpagelayout=OneColumn,%
pdfstartview=FitH,%
pdftitle={Deux extensions de Théorèmes de Hamburger},%
pdfauthor={Jean-François Burnol},%
pdfkeywords={Zeta function, Functional equation, Dirichlet series},%
pdfsubject={Number theory, Harmonic analysis}}

\usepackage{bookmark}
\usepackage{txfonts}

\renewcommand{\familydefault}{\ttdefault}
\renewcommand{\ttdefault}{txvtt}

\makeatletter

  \begingroup
  \nfss@catcodes
\DeclareFontFamily{T1}{txvtt}{\hyphenchar\font=127%
\fontdimen2\font=0.33333\fontdimen6\font%
\fontdimen3\font=0.16666\fontdimen6\font%
\fontdimen4\font=0.11111\fontdimen6\font}%

\DeclareFontShape{T1}{txvtt}{m}{n}{	
     <->t1xtt%
}{}%
\DeclareFontShape{T1}{txvtt}{m}{sc}{	
     <->t1xttsc%
}{}%
\DeclareFontShape{T1}{txvtt}{m}{sl}{	
     <->t1xttsl%
}{}%
\DeclareFontShape{T1}{txvtt}{m}{it}{	
     <->t1xttsl%
}{}%
\DeclareFontShape{T1}{txvtt}{m}{ui}{   	
     <->t1xttsl%
}{}%
\DeclareFontShape{T1}{txvtt}{bx}{n}{	
     <->t1xbtt%
}{}%
\DeclareFontShape{T1}{txvtt}{bx}{sc}{	
     <->t1xbttsc%
}{}%
\DeclareFontShape{T1}{txvtt}{bx}{sl}{	
     <->t1xbttsl%
}{}%
\DeclareFontShape{T1}{txvtt}{bx}{it}{	
     <->t1xbttsl%
}{}%
\DeclareFontShape{T1}{txvtt}{bx}{ui}{  	
     <->t1xbttsl%
}{}%
\DeclareFontShape{T1}{txvtt}{b}{n}{	
     <->t1xbtt%
}{}%
\DeclareFontShape{T1}{txvtt}{b}{sc}{	
     <->t1xbttsc%
}{}%
\DeclareFontShape{T1}{txvtt}{b}{sl}{	
     <->t1xbttsl%
}{}%
\DeclareFontShape{T1}{txvtt}{b}{it}{   	
     <->t1xbttsl%
}{}%
\DeclareFontShape{T1}{txvtt}{b}{ui}{   	
     <->t1xbttsl%
}{}%
  \endgroup

\re@DeclareMathSymbol{\Gamma}{\mathalpha}{lettersA}{0}
\let\alpha\alphaup
\let\beta\betaup
\let\delta\deltaup
\let\epsilon\epsilonup
\let\zeta\zetaup
\let\lambda\lambdaup
\let\mu\muup
\let\pi\piup
\let\sigma\sigmaup
\let\tau\tauup
\let\phi\phiup
\let\chi\chiup
\let\omega\omegaup

  \def\m@t@enc{\encodingdefault}
  \def\m@t@fam{\familydefault}
  \def\m@t@ser{\seriesdefault}
  \def\m@t@opsh{\shapedefault}  
  \def\m@t@bold{\bfdefault}
\DeclareSymbolFont{mtoperatorfont}
    {\m@t@enc}{\m@t@fam}{\m@t@ser}{\m@t@opsh}
\DeclareSymbolFontAlphabet{\mathrm}{mtoperatorfont}
\DeclareMathAlphabet{\mathbf}{\m@t@enc}{\m@t@fam}{\m@t@bold}{\m@t@opsh}
  \SetSymbolFont{mtoperatorfont}{bold} {\m@t@enc}
                                       {\m@t@fam}
                                       {\m@t@ser}
                                       {\m@t@opsh}
  \SetMathAlphabet{\mathbf}{bold}{\m@t@enc}
                                 {\m@t@fam}
                                 {\m@t@ser}
                                 {\m@t@opsh}
\def\operator@font{\mathgroup\symmtoperatorfont}

\DeclareMathSymbol{a}{\mathalpha}{mtoperatorfont}{`a}
\DeclareMathSymbol{b}{\mathalpha}{mtoperatorfont}{`b}
\DeclareMathSymbol{c}{\mathalpha}{mtoperatorfont}{`c}
\DeclareMathSymbol{d}{\mathalpha}{mtoperatorfont}{`d}
\DeclareMathSymbol{e}{\mathalpha}{mtoperatorfont}{`e}
\DeclareMathSymbol{f}{\mathalpha}{mtoperatorfont}{`f}
\DeclareMathSymbol{g}{\mathalpha}{mtoperatorfont}{`g}
\DeclareMathSymbol{h}{\mathalpha}{mtoperatorfont}{`h}
\DeclareMathSymbol{i}{\mathalpha}{mtoperatorfont}{`i}
\DeclareMathSymbol{j}{\mathalpha}{mtoperatorfont}{`j}
\DeclareMathSymbol{k}{\mathalpha}{mtoperatorfont}{`k}
\DeclareMathSymbol{l}{\mathalpha}{mtoperatorfont}{`l}
\DeclareMathSymbol{m}{\mathalpha}{mtoperatorfont}{`m}
\DeclareMathSymbol{n}{\mathalpha}{mtoperatorfont}{`n}
\DeclareMathSymbol{o}{\mathalpha}{mtoperatorfont}{`o}
\DeclareMathSymbol{p}{\mathalpha}{mtoperatorfont}{`p}
\DeclareMathSymbol{q}{\mathalpha}{mtoperatorfont}{`q}
\DeclareMathSymbol{r}{\mathalpha}{mtoperatorfont}{`r}
\DeclareMathSymbol{s}{\mathalpha}{mtoperatorfont}{`s}
\DeclareMathSymbol{t}{\mathalpha}{mtoperatorfont}{`t}
\DeclareMathSymbol{u}{\mathalpha}{mtoperatorfont}{`u}
\DeclareMathSymbol{v}{\mathalpha}{mtoperatorfont}{`v}
\DeclareMathSymbol{w}{\mathalpha}{mtoperatorfont}{`w}
\DeclareMathSymbol{x}{\mathalpha}{mtoperatorfont}{`x}
\DeclareMathSymbol{y}{\mathalpha}{mtoperatorfont}{`y}
\DeclareMathSymbol{z}{\mathalpha}{mtoperatorfont}{`z}
\DeclareMathSymbol{A}{\mathalpha}{mtoperatorfont}{`A}
\DeclareMathSymbol{B}{\mathalpha}{mtoperatorfont}{`B}
\DeclareMathSymbol{C}{\mathalpha}{mtoperatorfont}{`C}
\DeclareMathSymbol{D}{\mathalpha}{mtoperatorfont}{`D}
\DeclareMathSymbol{E}{\mathalpha}{mtoperatorfont}{`E}
\DeclareMathSymbol{F}{\mathalpha}{mtoperatorfont}{`F}
\DeclareMathSymbol{G}{\mathalpha}{mtoperatorfont}{`G}
\DeclareMathSymbol{H}{\mathalpha}{mtoperatorfont}{`H}
\DeclareMathSymbol{I}{\mathalpha}{mtoperatorfont}{`I}
\DeclareMathSymbol{J}{\mathalpha}{mtoperatorfont}{`J}
\DeclareMathSymbol{K}{\mathalpha}{mtoperatorfont}{`K}
\DeclareMathSymbol{L}{\mathalpha}{mtoperatorfont}{`L}
\DeclareMathSymbol{M}{\mathalpha}{mtoperatorfont}{`M}
\DeclareMathSymbol{N}{\mathalpha}{mtoperatorfont}{`N}
\DeclareMathSymbol{O}{\mathalpha}{mtoperatorfont}{`O}
\DeclareMathSymbol{P}{\mathalpha}{mtoperatorfont}{`P}
\DeclareMathSymbol{Q}{\mathalpha}{mtoperatorfont}{`Q}
\DeclareMathSymbol{R}{\mathalpha}{mtoperatorfont}{`R}
\DeclareMathSymbol{S}{\mathalpha}{mtoperatorfont}{`S}
\DeclareMathSymbol{T}{\mathalpha}{mtoperatorfont}{`T}
\DeclareMathSymbol{U}{\mathalpha}{mtoperatorfont}{`U}
\DeclareMathSymbol{V}{\mathalpha}{mtoperatorfont}{`V}
\DeclareMathSymbol{W}{\mathalpha}{mtoperatorfont}{`W}
\DeclareMathSymbol{X}{\mathalpha}{mtoperatorfont}{`X}
\DeclareMathSymbol{Y}{\mathalpha}{mtoperatorfont}{`Y}
\DeclareMathSymbol{Z}{\mathalpha}{mtoperatorfont}{`Z}

\DeclareMathSymbol{0}{\mathalpha}{mtoperatorfont}{`0}
\DeclareMathSymbol{1}{\mathalpha}{mtoperatorfont}{`1}
\DeclareMathSymbol{2}{\mathalpha}{mtoperatorfont}{`2}
\DeclareMathSymbol{3}{\mathalpha}{mtoperatorfont}{`3}
\DeclareMathSymbol{4}{\mathalpha}{mtoperatorfont}{`4}
\DeclareMathSymbol{5}{\mathalpha}{mtoperatorfont}{`5}
\DeclareMathSymbol{6}{\mathalpha}{mtoperatorfont}{`6}
\DeclareMathSymbol{7}{\mathalpha}{mtoperatorfont}{`7}
\DeclareMathSymbol{8}{\mathalpha}{mtoperatorfont}{`8}
\DeclareMathSymbol{9}{\mathalpha}{mtoperatorfont}{`9}

\DeclareMathSymbol{!}{\mathclose}{mtoperatorfont}{"21}
\DeclareMathSymbol{?}{\mathclose}{mtoperatorfont}{"3F}
\DeclareMathSymbol{*}{\mathalpha}{mtoperatorfont}{"2A}
\DeclareMathSymbol{,}{\mathpunct}{mtoperatorfont}{"2C}
\DeclareMathSymbol{.}{\mathord}{mtoperatorfont}{"2E}
\DeclareMathSymbol{:}{\mathrel}{mtoperatorfont}{"3A} 
\DeclareMathSymbol{;}{\mathpunct}{mtoperatorfont}{"3B}
\edef\mt@minus@sign{\mathord{\expandafter\mathchar\number\mathcode`\-}}
\def\relbar{\mathrel{\smash\mt@minus@sign}}
\def\rightarrowfill{$\m@th\mt@minus@sign\mkern-7mu  %
\cleaders\hbox{$\mkern-2mu\mt@minus@sign\mkern-2mu$}\hfill
  \mkern-7mu\mathord\rightarrow$}
\def\leftarrowfill{$\m@th\mathord\leftarrow\mkern-7mu %
 \cleaders\hbox{$\mkern-2mu\mt@minus@sign\mkern-2mu$}\hfill
  \mkern-7mu\smash\mt@minus@sign$}
\DeclareMathSymbol{-}{\mathbin}{mtoperatorfont}{21}
\DeclareMathSymbol{+}{\mathbin}{mtoperatorfont}{"2B}
\edef\mt@equal@sign{{\expandafter\mathchar\number\mathcode`\=}}
\DeclareRobustCommand\Relbar{\mathrel{\mt@equal@sign}}
\DeclareMathSymbol{=}{\mathrel}{mtoperatorfont}{"3D}
\DeclareMathDelimiter{(}{\mathopen} {mtoperatorfont}{"28}{largesymbols}{"00}
\DeclareMathDelimiter{)}{\mathclose}{mtoperatorfont}{"29}{largesymbols}{"01}
\DeclareMathDelimiter{[}{\mathopen} {mtoperatorfont}{"5B}{largesymbols}{"02}
\DeclareMathDelimiter{]}{\mathclose}{mtoperatorfont}{"5D}{largesymbols}{"03}
\DeclareMathDelimiter{/}{\mathord}{mtoperatorfont}{"2F}{largesymbols}{"0E}
\DeclareMathSymbol{/}{\mathord}{mtoperatorfont}{"2F}
\DeclareMathDelimiter{<}{\mathopen}{mtoperatorfont}{"3C}{largesymbols}{"0A}
\DeclareMathDelimiter{>}{\mathclose}{mtoperatorfont}{"3E}{largesymbols}{"0B}
\DeclareMathSymbol{<}{\mathrel}{mtoperatorfont}{"3C}
\DeclareMathSymbol{>}{\mathrel}{mtoperatorfont}{"3E}
\expandafter\DeclareMathDelimiter\@backslashchar
                        {\mathord}{mtoperatorfont}{"5C}{largesymbols}{"0F}
\DeclareMathDelimiter{\backslash}
    {\mathord}{mtoperatorfont}{"5C}{largesymbols}{"0F}
\DeclareMathSymbol\setminus\mathbin{mtoperatorfont}{"5C}
\DeclareMathSymbol{|}\mathord{mtoperatorfont}{"7C}
\DeclareMathDelimiter{|}{mtoperatorfont}{"7C}{largesymbols}{"0C}
\DeclareMathDelimiter\vert
                 \mathord{mtoperatorfont}{"7C}{largesymbols}{"0C}
\DeclareMathSymbol\mid\mathrel{mtoperatorfont}{"7C}
\DeclareMathDelimiter{\lbrace}
   {\mathopen}{mtoperatorfont}{"7B}{largesymbols}{"08}
\DeclareMathDelimiter{\rbrace}
   {\mathclose}{mtoperatorfont}{"7D}{largesymbols}{"09}

\DeclareMathSizes{9}{9}{7}{5}
\DeclareMathSizes{\@xpt}{\@xpt}{8}{6}
\DeclareMathSizes{\@xipt}{\@xipt}{9}{7}
\DeclareMathSizes{\@xiipt}{\@xiipt}{10}{8}
\DeclareMathSizes{\@xivpt}{\@xivpt}{\@xiipt}{10}
\DeclareMathSizes{\@xviipt}{\@xviipt}{\@xivpt}{\@xiipt}
\DeclareMathSizes{\@xxpt}{\@xxpt}{\@xviipt}{\@xivpt}
\DeclareMathSizes{\@xxvpt}{\@xxvpt}{\@xxpt}{\@xviipt}
\makeatother

\let\oldexists\exists
\renewcommand\exists{\oldexists\,}
\let\oldforall\forall
\renewcommand\forall{\oldforall\,}

\newtheorem{theo}{Théorème}

\newtheorem*{theo*}{Théorème}
\newtheorem*{cor1}{Corollaire 1}
\newtheorem*{cor2}{Corollaire 2}

\renewcommand{\Re}{\mathrm{Re}}
\renewcommand{\Im}{\mathrm{Im}}
\newcommand{\Oh}{\mathcal{O}}
\newcommand{\cF}{\mathcal{F}}
\newcommand{\NN}{\mathbf{N}}
\newcommand{\ZZ}{\mathbf{Z}}
\newcommand{\RR}{\mathbf{R}}

\newcommand{\CC}{\mathbf{C}}

\renewcommand\le\leqslant
\renewcommand\ge\geqslant

\newbox\boiteresume
\newenvironment{resume}{\global\setbox\boiteresume=\vbox\bgroup
  \hsize=\textwidth\def\baselinestretch{1}\small
  \noindent\unskip\textbf{Résumé}
 \par\medskip\noindent\unskip\ignorespaces}
 {\egroup}

\newbox\boiteabstract
\renewenvironment{abstract}{\global\setbox\boiteabstract=\vbox\bgroup
  \hsize=\textwidth\def\baselinestretch{1}\small
  \noindent\unskip\textbf{Abstract}
 \par\medskip\noindent\unskip\ignorespaces}
 {\egroup}

\newbox\boiteclefs
\newenvironment{clefs}{%
  \global\setbox\boiteclefs=\vbox\bgroup
  \hsize=\textwidth\def\baselinestretch{1}\small
  \normalfont\parskip0pt
  \ignorespaces}
{\egroup}

\newcommand*{\monheader}{{\normalfont
\ifodd\value{page}%
   \emph{sur l'équation fonctionnelle de la fonction $\zeta$ de Riemann}
  \hfil\thepage
\else \thepage
  \hfil\emph{Deux extensions de Théorèmes de Hamburger}\fi}}

\begin{document}
\pdfbookmark[1]{Page de titre}{pagedetitre}

\thispagestyle{empty}

\markright{\protect\monheader}
\makeatletter
\def\@oddhead{\rightmark}
\makeatother

\begin{flushright} 
\small\ladate\hspace{\the\leftmargin}
\end{flushright}

\long\def\letitre{Deux extensions de Théorèmes de Hamburger\\
\centerline{\hss\mbox{(portant sur l'équation fonctionnelle de
  la fonction zêta)}\hss}}

\long\def\letitreanglais{%
Two extensions of Hamburger's theorems\\
(on the functional equation of the zeta function)}

\def\moi{Jean-François~Burnol}

\def\merci{L'auteur remercie le C.R.M. de Barcelone pour l'hospitalité de son
  accueil en mai 2011, lors d'un séjour pendant lequel ce travail a été conçu
  et rédigé.}

\def\courriel{burnol@math.univ-lille1.fr}

\def\monadresse{Université Lille 1,  
UFR de Mathématiques,
Cité Scientifique M2,
F-59655 Villeneuve d'Ascq,
France}

\begin{resume}
  Nous proposons deux types d'extensions aux théorèmes de
  Hamburger sur les séries de Dirichlet avec équation
  fonctionnelle comme celle de la fonction zêta de Riemann,
  sous des hypothèses plus faibles. Ceci repose sur le
  dictionnaire entre les fonctions méromorphes modérées avec
  cette équation fonctionnelle et les distributions
  tempérées avec la condition de support
  $S$-étendue.
\end{resume}

\begin{abstract}
  We propose two types of extensions to Hamburger's theorems on the Dirichlet
  series with functional equation like the one of the Riemann zeta function,
  under weaker hypotheses. This builds upon the dictionary between the
  moderate meromorphic functions with functional equation and the tempered
  distributions with extended $S$-support condition.
\end{abstract}

\begin{clefs}
  \noindent\textit{Keywords\string:} Riemann Zeta function, Dirichlet
  series, Hamburger Theorem, functional equations, co-Poisson formula, Poisson
  formula.
  \par
  \noindent\textit{MSC2000\string:} 11M06, 11F66.
\end{clefs}

\begin{center}
  \def\baselinestretch{1}%
    \Large\letitre\par\vskip18pt
    \normalsize\moi\footnote{\texttt{\courriel}}\footnote{\merci}\par\vskip10pt
    \footnotesize\itshape\centerline{\hss\mbox{\monadresse}\hss}\vskip18pt
    \hrule\vskip12pt
\pdfbookmark[2]{Résumé}{resumei}%
    \unvbox\boiteresume\par\vskip16pt
\pdfbookmark[2]{Title and abstract}{resumeii}%
    \large\normalfont\letitreanglais\par\vskip10pt
    \unvbox\boiteabstract\par\vskip10pt
    \unvbox\boiteclefs\par\vskip10pt
    \hrule\vskip24pt
\end{center}

\section{Introduction et présentation des résultats}

Le Théorème de Hamburger \cite{hamb1} dit à peu près que deux
fonctions $f(s)$ et $g(s)$, méromorphes dans le plan complexe, qui
admettent chacune pour $\Re(s)\gg1$ une représentation sous 
forme de série de Dirichlet $\sum_{n=1}^\infty a_n n^{-s}$, et
sont reliées par l'équation fonctionnelle
 \[ {\pi^{-\frac s2}\Gamma(\frac s2)} g(s) =
{\pi^{-\frac{1-s}2}\Gamma(\frac{1-s}2)} f(1-s)\;, \] sont alors
nécessairement égales à un multiple de la fonction zêta de Riemann
$\zeta(s)$. Ce résultat fameux illustre une certaine rigidité de
l'équation fonctionnelle de la fonction zêta. Notre but ici est,
en utilisant quelques notions communes de la théorie des
distributions, de le reprouver sous des hypothèses nettement
plus faibles que celles d'origine. En particulier nous
n'aurons pas besoin de demander à $g$ d'admettre une représentation
en série de Dirichlet, mais seulement de tendre rapidement vers
une limite lorsque $\Re(s)\to\infty$. Et lorsque l'on supposera 
seulement $g$ bornée dans un demi-plan alors $f$ sera nécessairement
une combinaison linéaire finie de
$\zeta(s)$, $\zeta(s+2)$, $\zeta(s+4)$, \dots{}~.

L'article se veut accessible à tout lecteur disposant du bagage
usuel de base de la théorie des distributions: il faut y ajouter
quelques éléments qui ont été développés dans le chapitre IV de
\cite{jfcopoisson}, chapitre qui peut être lu avec les mêmes
pré-requis: il y est construit une notion de «fonction méromorphe
modérée avec équation fonctionnelle» dont nous rappellerons les
principaux éléments.

Il est plus commode pour cet article de mettre l'équation
fonctionnelle de la fonction zêta sous la forme $\zeta(s) =
\chi(s) \zeta(1-s)$ avec
\[ \chi(s) = \frac{\pi^{-\frac{1-s}2}\Gamma(\frac{1-s}2)}{\pi^{-\frac
    s2}\Gamma(\frac s2)}\] La fonction méromorphe $\chi(s)$ est à croissance
au plus polynomiale dans toute bande verticale de largeur
finie. Voici tout d'abord l'énoncé originel démontré par Hamburger:

\pdfbookmark[2]{Premier théorème de Hamburger}{hambi}

\begin{theo*}[Hamburger, \cite{hamb1}]
  Soit $f$ une fonction méromorphe dans le plan complexe tout entier, et
  d'ordre fini ($f$ est $\Oh(\cramped{e^{|s|^k}})$ avec un certain entier $k$
  pour $|s|\to\nobreak\infty$, en particulier ne possède au plus qu'un nombre
  fini de pôles). Si $f(s)$ est représentée pour $\Re(s)>\nobreak1$ par une
  série de Dirichlet $\sum_{n=1}^\infty {a_n}{n^{-s}}$ absolument
  convergente, et si la fonction méromorphe \[g(s) = \chi(s)f(1-s)\] admet
  elle aussi, pour $\Re(s)\gg1$, une représentation sous la
  forme d'une série de Dirichlet convergente $\sum_{n=1}^\infty
  {b_n}{n^{-s}}$, alors $f$ est un multiple de la fonction zêta. En
  particulier si $f(s)$ est représentée pour $\Re(s)>\nobreak1$ par une série
  de Dirichlet absolument convergente $\sum_{n=1}^\infty {a_n}{n^{-s}}$ et
  vérifie l'équation fonctionnelle
 \[
\pi^{-\frac
    s2}\Gamma(\frac s2)f(s) = \pi^{-\frac{1-s}2}\Gamma(\frac{1-s}2) f(1-s)\;,
\]
alors elle est un multiple de la fonction zêta de Riemann. 
\end{theo*}

Énonçons maintenant notre résultat principal:
\begin{theo}\label{th:1}
\addcontentsline{toc}{subsection}{\autoref{th:1}}
Soit $f(s) = \sum_{n=1}^\infty \frac{a_n}{n^s}$ une série de Dirichlet qui
converge pour $\Re(s)$ suffisamment grand. On suppose que:
  \begin{enumerate}[\upshape(1)]
  \item\label{cond:1} elle admet un prolongement méromorphe au plan complexe
    tout entier, avec un nombre fini de pôles,
  \item\label{cond:2} dans toute bande verticale de largeur finie et pour tout $\epsilon>0$
    on a $f(s) = \Oh(e^{\exp{\epsilon |s|}})$ lorsque $|\Im(s)|\to\infty$,
  \item\label{cond:3} la fonction méromorphe $g(s) = \chi(s) f(1-s)$ vérifie
\[ \exists N\in\NN, \exists y>\frac12, \qquad g(s) = \Oh(|s|^N y^{-s})\qquad
\text{pour }\Re(s)\to+\infty\;.\]
  \end{enumerate}
Alors $f(s)$ est une somme finie:
 \[f(s) = \sum_{k\in\ZZ} c_k \zeta(s-2k)\]
\end{theo}

Nous voyons que le \autoref{th:1} fait des hypothèses analytiques
nettement plus faibles en ce qui concerne le comportement
analytique de $f(s)$ et celui de $\chi(s)f(1-\nobreak s)$. La
conclusion n'est plus, dans l'immédiat,
que $f$ est zêta, mais il suffit d'être un peu plus exigeant envers la
fonction $g$:

\pdfbookmark[2]{Corollaires du \autoref{th:1}}{cor}%
\begin{cor1}
  Avec les notations du \emph{\autoref{th:1}}, si \textup{(\ref{cond:1})},
  \textup{(\ref{cond:2})}, et si
  \begin{enumerate}[\upshape(1${}'$)] \setcounter{enumi}2
  \item la fonction méromorphe $g(s) = \chi(s) f(1-s)$ est bornée pour
    $\Re(s)\gg 1$ ou même vérifie seulement 
\[ \exists y>\frac12\qquad g(s) \mathop{=}_{\Re(s)\to+\infty} \Oh(y^{-s})\;,\]
 \end{enumerate}
alors $f$ est une combinaison linéaire finie de $\zeta(s)$, $\zeta(s+2)$,
 ${\zeta(s+4)}$~\dots{}
\end{cor1}
\begin{cor2}
 Si \textup{(\ref{cond:1})}, \textup{(\ref{cond:2})}, \textup{(\ref{cond:3})},
 et
  \begin{enumerate}[\upshape(1)] \setcounter{enumi}3
  \item\label{cond:4} il existe $c$ tel que pour sigma réel tendant
    vers $+\infty$ on a $g(\sigma) = c + \Oh(\sigma^{-k})$ pour tout
    $k\in\NN$,
 \end{enumerate}
alors $f(s) = c\, \zeta(s)$.
\end{cor2}

Un deuxième théorème de Hamburger englobe son premier. Il demande toujours que
$f(s)$ soit une série de Dirichlet au sens strict, mais autorise $g(s)$ à être
une série de Dirichlet générale, c'est-à-dire de la forme
$\sum_{n=1}^\infty \frac{b_n}{\lambda_n^s}$ avec $0<\lambda_n\to+\infty$.
\pdfbookmark[2]{Second théorème de Hamburger}{hambii}
\begin{theo*}[Hamburger, \cite{hamb2}]
  Soit $f$ une fonction méromorphe dans le plan complexe tout entier, d'ordre
  fini, et égale pour $\Re(s)>1$ à la somme d'une série de Dirichlet
  absolument convergente. Si la fonction \[g(s) = \chi(s)f(1-s)\] admet pour
  $\Re(s)$ suffisamment grand une représentation sous la forme d'une série de
  Dirichlet générale\footnote{On pourra toujours supposer que les entiers font
    partie des $y_n$; mais on n'autorisera $b_n=0$ \emph{que} lorsque
    $y_n\in \NN$.}
  \[ g(s) = \sum_{n=1}^\infty \frac{b_n}{y_n^s},\quad \Re(s)\gg0,
  0<y_n\to\infty, 0< y_1 < y_2 < \dots < y_k = 1 < \dots \] alors $y_{n+k} =
  y_n+1$ pour tout $n\ge1$, avec $b_{n+k} = b_n$ et de plus, pour $1\le j < k$
  on a les symétries $b_j = b_{k-j}$, $y_j+y_{k-j} = 1$. Ainsi, si tous les
  $y_n$ sont $\ge1$, c'est que $k=1$ et que $f$ est un multiple de la fonction
  zêta. En général on a
\[ f(s) = \sum_{n=1}^\infty \frac{\sum_{1\le j \le k} b_j\cos(2\pi y_j
  n)}{n^s},\qquad 0<y_1<\dots <y_k = 1\] 
\end{theo*}

Il subsiste dans cet énoncé une certaine dissymétrie entre la
série $f(s)$ formée avec des entiers et la série $g(s)$. Nous
ajouterons donc dans le présent article à notre Théorème principal
le suivant qui rétablit la symétrie entre $f$ et $g$:

\begin{theo}\label{th:2}
  \addcontentsline{toc}{subsection}{\autoref{th:2}} 
  Soit \[f(s) = \sum_{n=1}^\infty \frac{a_n}{x_n^s}\] une série de Dirichlet
  générale ($0 < x_n < \infty$, $x_n\to\infty$) qui converge pour $\Re(s)$
  suffisamment grand. On suppose que:
  \begin{enumerate}[\upshape(1)]
  \item $f$ admet un prolongement méromorphe au plan complexe tout entier,
    avec un nombre fini de pôles,
  \item dans toute bande verticale de largeur finie on a
    $f(s) = \Oh(e^{\exp{\epsilon |s|}})$ pour tout $\epsilon>0$ lorsque
    $|\Im(s)|\to\infty$,
  \item la fonction méromorphe $g(s) = \chi(s) f(1-s)$ est aussi, pour
    $\Re(s)\gg0$ la somme d'une série de Dirichlet générale $\sum_{n=1}^\infty
    \frac{b_n}{y_n^s}$.
  \end{enumerate}
Alors les assertions suivantes sont équivalentes:
\begin{enumerate}[\upshape(i)]
\item l'ensemble $\{x_n | a_n\neq 0\}$ est contenu dans un nombre fini de
  progressions arithmétiques de raison $1$,
\item l'ensemble $\{y_n | b_n\neq 0\}$ est contenu dans un nombre fini de
  progressions arithmétiques de raison $1$,
\item il existe un nombre fini de couples réels $(d,e)$ tels que la série
  $f(s)$ soit de la forme:
\[ f(s) = \sum  c(d,e) \Biggl(\smashoperator[r]{\sum_{\substack{x\equiv
    d\bmod1\\x >0}}} \frac{\exp(2\pi i ex)}{x^s} + \smashoperator[r]{\sum_{\substack{x\equiv
    -d\bmod1\\x >0}}} \frac{\exp(-2\pi i ex)}{x^s}\Biggr)\]
\end{enumerate}
\end{theo}

Si l'on connaît a priori le nombre $N$ de suites arithmétiques nécessaires pour
$f$ alors on n'a plus besoin de supposer que $g$ soit aussi une série
de Dirichlet générale, il suffit d'en contrôler le \og début\fg:

\begin{theo}\label{th:3}
  \addcontentsline{toc}{subsection}{\autoref{th:3}}
  Soit \[f(s) = \sum_{n=1}^\infty \frac{a_n}{x_n^s}\] une série de Dirichlet
  générale ($0 < x_n < \infty$, $x_n\to\infty$) qui converge pour $\Re(s)$
  suffisamment grand. On suppose que:
  \begin{enumerate}[\upshape(1)]
  \item $f$ admet un prolongement méromorphe au plan complexe tout entier,
    avec un nombre fini de pôles,
  \item dans toute bande verticale de largeur finie on a
    $f(s) = \Oh(e^{\exp{\epsilon |s|}})$ pour tout $\epsilon>0$ lorsque
    $|\Im(s)|\to\infty$,
  \item l'ensemble $\{\pm x_n | a_n\neq0\}\subset \RR$ est contenu dans un
    nombre fini 
    $N$ de translatés de $\ZZ$,
\item il existe $Y>N$, ainsi que $(y_n,b_n)\in
  \OO{0,+\infty}\times\CC$, $1\le n \le M$, tels que
\[ g(s) \coloneqq \chi(s) f(1-s) \mathop{=}\limits_{\Re(s)\to+\infty}
\sum_{1\le n\le M} b_n y_n^{-s} + \Oh(Y^{-s})\;.\] 
  \end{enumerate}
  Alors en fait $g$ est la somme d'une série de Dirichlet générale et
  conformément au théorème précédent il existe
  une représentation:
\[ f(s) = \sum c(d,e) \Biggl(\smashoperator[r]{\sum_{\substack{x\equiv
    d\bmod1\\x >0}}} \frac{\exp(2\pi i ex)}{x^s} + \smashoperator[r]{\sum_{\substack{x\equiv
    -d\bmod1\\x >0}}} \frac{\exp(-2\pi i ex)}{x^s}\Biggr)\]
\end{theo}

\pdfbookmark[2]{Discussion}{discussion}

Voici quelques références, par nécessité très brèves, à la
littérature.  Les résultats de Hamburger \cite{hamb1,hamb2,hamb3}
(qui incluent aussi des énoncés semblables avec des équations
fonctionnelles formées avec
\smash{$\pi^{-\frac{s+1}2}\Gamma(\frac{s+1}2)$}) sont bien connus
et ont eu une postérité riche et assez variée. En particulier les
travaux de Hecke \cite{hecke} et de Weil \cite{weil}, sur le lien
important avec l'invariance modulaire sur le demi-plan de
Poincaré, et donc d'une manière générale avec les fonctions
automorphes, ont eu une forte influence motivant l'établissement
de «converse theorems» de plus en plus sophistiqués dans le cadre
du programme de Langlands. Également et de manière liée, l'étude
des axiomes de Selberg est intimement concernée par des théorèmes
inverses comme celui de Hamburger: ainsi Kaczorowski, Molteni et
Perelli ont caractérisé les conducteurs pour lesquels les
fonctions de Dirichlet $L(s,\chi)$ sont les seules avec leur
équation fonctionnelle dans une certaine classe \cite{kmp}. Il
s'agit là de questions d'une nature infiniment plus arithmétique
que celles que nous examinerons ici-même, puisqu'il faut
déterminer les caractères de Dirichlet partageant la même parité
et la même somme de Gauss. Molteni a obtenu dans \cite{molt1,
  molt2} des résultats étendant ceux de \cite{kmp} à une classe
plus large de conducteurs.

Pour en revenir au contexte originel, Hamburger avait déjà insisté
\cite{hamb4} sur la variété des écritures équivalentes de
l'équation fonctionnelle comme formule sommatoire.  On peut y voir
comme une anticipation du point de vue des distributions, qui
déplace la focalisation de la \og fonction test\fg{} vers la forme
linéaire sur ces fonctions tests (même si on doit garder présent à
l'esprit que ce n'est là qu'une façon d'appréhender la réalité
d'une distribution). Ceci est dit uniquement dans la perspective
de l'équation fonctionnelle, car pour d'autres propriétés, on ne
saurait préjuger de ce qui s'avérera plus ou moins utile. Par
exemple chez Hamburger \cite{hamb4} et Siegel \cite{siegel}
l'équation fonctionnelle devient un développement en éléments
simples d'une certaine transformée de Laplace.

Cette écriture fut généralisée par Bochner et Chandrasekharan
\cite{bc}; et par Chandrasekahran et Mandelbrojt \cite{cm1,
  cm2}. On considère des séries de Dirichlet $\sum a_n
\lambda_n^{-s}$ et $\sum b_n \mu_n^{-s}$ reliées par des équations
fonctionnelles d'un certain type et on pose la question: étant
donnés les $\lambda_n$ et les $\mu_n$ peut-on en fonction de leurs
propriétés majorer la dimension de l'espace des choix possibles
pour les $a_n$ et $b_n$?  Kahane et Mandelbrojt (\cite{km})
interprètent l'équation fonctionnelle comme une formule de
transformée de Fourier pour des distributions \og
quasi-périodiques\fg. Ils prennent ensuite comme point de départ
une telle paire de Fourier, et prouvent dans ce contexte des
énoncés sur le support et le spectre: pour certains de leurs
résultats, ils reviennent à des fonctions analytiques, via une
transformée de Laplace-Carleman, pour d'autres ils mettent en plus
à l'oeuvre des techniques de fonctions quasi-périodiques. Sous
l'hypothèse d'une densité supérieure de répartition finie pour les
$\mu_n$, ils prouvent que les $\lambda_n$ sont combinaisons
linéaires à coefficients entiers d'un nombre fini de nombres réels
(ce qui généralise un résultat de \cite{cm1}). Ils ne s'attardent
pas à considérer des supports limités d'emblée à un nombre fini de
suites arithmétiques, cela ne serait chez eux qu'un exercice
facile et immédiat une fois dans le contexte des distributions,
mais obtiennent des critères qui a posteriori donnent cette
situation.\footnote{Chez nous, les suites arithmétiques sont
  restreintes à la raison $1$.}
Lorsque les $\lambda_n$ sont a priori des entiers (ou placés sur un nombre
fini de translatés de $\ZZ$) la difficulté est  de faire
la traduction vers les distributions: car après il n'est même plus alors besoin
d'outils avancés. 

Après tous les travaux cités, et d'autres encore certainement, il
est donc devenu de plus en plus «bien connu» que l'équation
fonctionnelle traduit l'invariance de la distribution
$\sum_{n\in\ZZ} \delta(x-n)$ sous la transformation de Fourier
\footnote{La transformation de Fourier est définie ici avec le
  noyau $\exp(2\pi i xy)$. Suivant que l'on étudie l'équation
  fonctionnelle formée avec $\Gamma(\frac s2)$ ou
  $\Gamma(\frac{s+1}2)$ il faudra considérer des fonctions et
  distributions paires ou impaires.}. Mais c'est avec l'idée de
co-Poisson \cite{jffourierzeta} que nous obtenons une vision plus
complète de cette correspondance. C'est sur la base de l'idée de
co-Poisson que le chapitre IV de \cite{jfcopoisson} établit un
dictionnaire (allant dans les deux sens) entre une certaine classe
de fonctions méromorphes et une certaine classe de distributions
tempérées.


Des questions fines (comme le résultat de Kahane-Mandelbrojt) apparaissent
lorsque l'on veut étudier en général des mesures ou des distributions à
support discret $D$ et dont la transformée de Fourier $\cF(D)$ est encore à
support discret\footnote{\emph{Cf.} \cite{cord} pour un résultat d'unicité
  dans $\RR^n$, lorsque l'on a deux mesures discrètes positives dont l'une est
  de la forme $\sum_k \delta_{x_k}$.}; l'objet du présent travail n'est
aucunement une telle étude générale. Au contraire, nous bénéficierons (comme
Hamburger) dès notre point de départ de la simplification énorme apportée par
l'hypothèse que la distribution $D = \sum_{n\ge1} a_{n} (\delta_{x_n} +
\delta_{-x_n})$ associée à la fonction analytique $f$ est une mesure supportée
sur les entiers $\ZZ$ ou un nombre fini de ses translatés. Mais le
dictionnaire général de \cite[Chap. IV]{jfcopoisson} nous donne les moyens
d'affaiblir considérablement les hypothèses faites par Hamburger et d'en
étendre ses énoncés, et c'est là tout l'objet de la présente rédaction.

\section{Preuves}

\pdfbookmark[2]{Un dictionnaire et la preuve du \autoref{th:1}}{preuvei}

Le chapitre IV de \cite{jfcopoisson} établit une correspondance générale via
la transformation de Mellin \emph{droite} \og $\int_0^\infty D(x)
x^{-s}\,dx$\fg{} entre des fonctions $f(s)$ avec un nombre fini de pôles dans
tout le plan complexe et une certaine condition de croissance (pour tout
demi-plan $\Re(s)>\sigma$ on a l'existence d'un entier $N$ et d'un réel $A>0$
avec, loin des pôles, $f(s) = O(|s|^N A^{\Re(s)})$, condition qui doit aussi
être imposée à $\chi(s) f(1-s)$), et, d'autre part, les distributions
tempérées paires $D(x)$, nulles dans un voisinage de l'origine et dont la
transformée de Fourier $E(x)$ est, non pas forcément nulle (cela serait le cas
si $f(s)$ et $\chi(s)f(1-s)$ n'avaient pas de pôles), mais
\emph{quasi-homogène} dans un voisinage de l'origine.

Une distribution $Q(x)$ est dite quasi-homogène si les $Q(tx)$, $t\neq 0$
engendrent un espace vectoriel de dimension finie. Par exemple $\log|x|$ est
quasi-homogène. Les distributions  \emph{homogènes} sont essentiellement les
$|x|^{w-1}$ (je ne m'occuperai que de distributions \emph{paires}), avec $w\in
C$, à ceci près que pour $w= 0$, $w=-2$, $w= -4$, \dots{}
il faut en fait prendre $\delta(x)$, $\delta''(x)$,
$\delta^{(4)}(x)$, \dots. Plus précisément, soit:
\[ Q_w(x) = \frac{|x|^{w-1}}{\pi^{-\frac w2}\Gamma(\frac w2)}\] Lorsque
$0<\Re(w)$ cela définit directement une distribution homogène de $x$, et pour
$\Re(w)\le 0$ on peut prouver que $Q_w$ se prolonge comme fonction analytique
à valeurs dans les distributions (tempérées). On notera que la restriction de
la distribution $Q_w$ à $x\neq 0$ est précisément la fonction donnée par la
formule ci-dessus, donc identiquement zéro si $\frac w2$ est un pôle de la
fonction Gamma.\footnote{La fonction $1/|x|$ n'est pas une distribution mais
  il y a une infinité de distributions qui se restreignent sur $x\neq0$ à
  cette fonction; parmi celles-ci il y a la transformée de Fourier de
  $-2\log|x|$.
Ces distributions sont toutes
  quasi-homogènes, mais aucune n'est homogène, et la seule (à un multiple
  près) distribution avec l'homogénéité de $1/|x|$ c'est le Dirac à
  l'origine. Par contre la fonction $1/x$ elle est une distribution homogène;
  elle est impaire, et c'est la valeur principale de Cauchy.} On dispose de la
formule utile $\cF(Q_w) = Q_{1-w}$, et par exemple $\cF(Q_0) = Q_1 = 1$ ce qui
permet de voir $Q_0 = \delta$.

Les homogènes sont des vecteurs propres de la dérivation $\frac{d}{dx} x$, et
on peut aussi définir les quasi-homogènes comme celles qui sont annulées par
des polynômes en cet opérateur. Pour plus de détails, voir
\cite[IV]{jfcopoisson}. 

Rappelons aussi ce que l'on entend par $\int_0^\infty D(x) x^{-s}\,dx$,
lorsque $D$ est une distribution tempérée paire. Nous ne faisons cette
définition que lorsque $D$ restreinte à un intervalle $\OO{-a,a}$ est
quasi-homogène.\footnote{Dans la vraie vie on a bien sûr besoin de calculer
  des transformées de Mellin plus générales, par exemple
  \smash{$\int_0^\infty f(x) x^{-s}\,dx$} lorsque $f\in L^2(0,+\infty;dx)$
  mais alors ce Mellin n'est ni moins ni mieux qu'une fonction de carré
  intégrable pour $\Re(s) = \frac12$. Si l'on veut de l'analyticité, il se
  trouve que la condition de support considérée est la plus simple ayant déjà
  une utilité. Le fait que la condition de support puisse faire bon ménage
  avec Fourier est déjà en soi un fait remarquable. }  Il est important que
les identités d'Euler $\sum_{n\in\ZZ} q^n = 0$ et $\int_0^\infty x^{s-1}dx =
0$ indiquent que les transformées de Mellin des distributions
quasi-homogènes sont identiquement nulles.  On commence donc par soustraire la
partie quasi-homogène de $D$ (qui pour $x>0$ est une fonction analytique
combinaison linéaire d'expressions $(\log x)^N x^{w-1}$, $N\in\NN$, $w\in\CC\setminus(-2\NN)$;
on fait la soustraction de cette expression analytique sur \emph{tout}
l'intervalle $\OO{0,+\infty}$), donc on peut supposer $D$ nulle sur $\OO{0,a}$ (il
peut rester des Dirac à l'origine, que l'on oublie). On peut alors écrire $D$
pour $x>0$ comme la dérivée $N\ieme$ d'une fonction continue $C$ à croissance
polynomiale, nulle sur $\OO{0,a}$, et on pose (\cite[4.11]{jfcopoisson}):
\[ \widehat D(s) = s(s+1)\dots(s+N-1) \int_a^\infty C(x) x^{-s-N}\,dx,\qquad
\Re(s)\gg0\] Toutes les formules raisonnables marchent avec cette
définition. Par exemple (\cite[4.20]{jfcopoisson}) si $D$ est de la forme
$\sum_n a_n (\delta(x - x_n) + \delta(x+x_n))$, $0<\nobreak x_1 < x_2 < \dots \to
\infty$, $\sum_{x_n \le X} |a_n| =\Oh( X^A)$, $A<\infty$ alors $\widehat D(s)
= \sum_{n\ge1} a_n x_n^{-s}$.

\bigskip

Venons-en à la démonstration du \autoref{th:1}.
Soit $D$ la distribution paire $\sum_{n\in Z} a_n \delta_n$ avec $a_0 = 0$ et
$a_{-n} = a_n$. Comme la série de Dirichlet converge quelque part les $a_n$
ont croissance polynomiale et $D$ est une distribution tempérée. Et nous
venons de préciser que la série de Dirichlet est aussi la transformée de
Mellin, en notre sens, de $D$.

Ensuite, par hypothèse $\chi(s) f(1-s)$ est pour $\Re(s) = \sigma$
 à croissance polynomiale lorsque $\sigma\gg0$ donc $f(s) =
\chi(s)\chi({1-s})f(s)$ est à croissance polynomiale sur toute droite
$\Re(s) = -\sigma$ avec $\sigma\gg0$. Si on va
suffisamment loin à droite $f(s)$ est bornée. Fixons donc deux telles droites
verticales $\Re(s) = -\sigma$ et $\Re(s) = +\sigma$ avec $\sigma\gg0$. Entre les deux on a 
l'hypothèse a priori (modulo un nombre fini de pôles) d'une croissance en
$\Oh(\exp( e^{\epsilon |s|}))$. Soit $P(s)$ un polynôme correspondant aux
pôles, prenons $A\gg0$ et $M$ entier suffisamment grand; le produit
$P(s)f(s)/(s+A)^M$ tend vers zéro sur chacune des deux droites $\Re(s) = \pm
\sigma$ et est $\Oh(\exp( e^{\epsilon |s|}))$ entre, pour tout
$\epsilon>0$. Par Phragmén-Lindelöf $P(s)f(s)/(s+A)^M$ est borné dans cette bande, et
$f$ y est à croissance polynomiale. Comme $f$ est donnée par une série de
Dirichlet elle est bornée pour $\Re(s) \gg 0$, donc $f$ est à croissance
polynomiale dans tout demi-plan droit. 

La fonction $f(s)$ est donc une fonction méromorphe \emph{modérée} au sens de
\cite[Déf. 4.27]{jfcopoisson}. Considérons la fonction $g(s) = \chi(s)
f(1-s)$. Par ce qui a déjà été établi pour $f$ elle est à croissance
polynomiale dans toute bande verticale de largeur finie. Et par hypothèse elle
est $\Oh(|s|^N Y^{-s})$, pour $\Re(s)\ge \sigma_0$, lorsque $\sigma_0\gg 0$
est bien choisi, avec $N$ entier et $Y>\frac12$. Les pôles de \smash{$g(s) =
  \pi^{s-\frac12} \Gamma(\frac{1-s}2)\Gamma(\frac s2)^{-1} f(1-s)$}, mis à
part ceux qui proviennent de $f$, ne peuvent être qu'en $1-s = 0$, $-2$, $-4$,
\dots, mais comme il n'y en a pas pour $\Re(s)\gg 0$, c'est donc que $g$ n'a
qu'un nombre fini de pôles dans tout le plan complexe. Soit $Q(s)$ un polynôme
correspondant à ces pôles, pour $\Re(s)\ge-\sigma_0$ la fonction analytique
$Y^{s} Q(s)g(s)$ est à croissance polynomiale: car elle l'est d'une part pour
$|\Re(s)|\le \sigma_0$ et d'autre part pour $\Re(s)>\sigma_0$. Et cela est
vrai pour tout $\sigma_0$, donc $g$ est aussi une fonction \emph{modérée} au
sens de \cite{jfcopoisson}.

Par \cite[Thm.{} 4.54]{jfcopoisson}: \emph{la distribution $D$ dont $f$ est la
  transformée de Mellin a une transformation de Fourier $E$ dont la
  restriction a un certain intervalle $\OO{-Y_0,Y_0}$ est quasi-homogène; de
  plus la transformée de Mellin de la distribution tempérée paire $E$ est
  égale à $g(s) = \chi(s) f(1-s)$.}

On dispose d'une information complémentaire importante \cite[Rem.{}
4.29]{jfcopoisson}: notons $Q$ la composante quasi-homogène de $E$, et soit
$E_1 = E - Q$, donc cet $E_1$ est identiquement nulle dans un certain
intervalle $\OO{-Y_0,Y_0}$. Si on prend pour $Y_0$ le plus petit point positif du
support de $E_1$ alors:
\[ -\log Y_0 = \limsup_{\sigma\to+\infty} \frac1{\sigma}\log |g(\sigma)|\]
Nous avons dans notre cas {\boldmath$Y_0 \ge Y >\frac12$}.

Comme $D$ est une mesure supportée sur $\ZZ$ elle vérifie l'équation 
\[ \exp(2\pi i x)D(x) = D(x)\]
Donc, et cela sera fondamental évidemment:
\[ E(x+1) = E(x)\quad\text{et ainsi } Q(x+1) - Q(x) = E_1(x) - E_1(x+1)\]

La restriction de la partie quasi-homogène $Q$ à $x\neq 0$ est une certaine
fonction $q(x)$. Cette fonction est paire. Pour $x>0$ elle est une expression
finie de la forme $q(x) = \sum_{n,w} q_{n,w} (\log x)^n x^{w-1}$, en
particulier elle est analytique sur $\CC\setminus\OF{-\infty,0}$ avec croissance
polynomiale à l'infini. Regardons $\epsilon$ avec $0<\epsilon<\frac12$,
$\frac12+\epsilon\le Y_0$. Sur $\frac12-\epsilon<x<\frac12+\epsilon$, $E_1(x)$
comme $E_1(x-1)$ sont identiquement nuls, donc
\[ \frac12-\epsilon<x<\frac12+\epsilon\implies q(x) = q(x-1) = q(1-x)\] La
fonction $q$ est analytique sur $\CC\setminus \OF{\infty,0}$ et la
fonction $q(1-x)$ est analytique sur $\CC\setminus
\FO{1,+\infty}$. Ces deux fonctions coïncident sur l'intervalle
$\OO{\frac12-\epsilon,\frac12+\epsilon}$, \emph{donc $q$ est une fonction
  entière.} Or elle est à croissance polynomiale. Donc $q$ est un polynôme sur
$x>0$. La conclusion est que $q(x) = P(|x|)$ avec $P$ un polynôme qui vérifie
$P(x) = P(1-x)$.

Considérons la fonction (continue) $1$-périodique $B(x)$ qui vaut $P(x)$ pour
$0<x<1$. Nous avons $B(x) = q(x)$ pour $0<x<1$ puis, pour $-1<x<0$: $B(x) =
B(x+1) = q(x+1) = q(x)$. Les polynômes de Bernoulli d'indices pairs $B_0 = 1$,
$B_2 = x^2 - x + \frac16$, \dots{} forment une base de l'espace vectoriel des
polynômes invariants par $x\mapsto 1-x$, et (à une constante multiplicative
près) la fonction $1$-périodique $B_{2n}(\{x\})$ (pour $n\ge1$) est la
transformée de Fourier de la distribution
$\sum_{k\in\ZZ\setminus\{0\}} k^{-2n}\delta_k$. Ainsi, quitte à retirer de
notre série de Dirichlet d'origine une combinaison linéaire finie de
$\zeta(s+\nobreak2)$, $\zeta(s+4)$, \dots{} on peut faire en sorte que la
fonction $1$-\nobreak périodique $B(x)$ soit réduite à une constante. La
transformée de Fourier de $\sum_{k\neq 0} \delta_k$ est $\sum_{k\neq 0}
\delta_k + \delta - 1$, et quitte à retirer en plus un multiple de $\zeta(s)$
(à la serie de Dirichlet originellement considérée) on peut réduire $B(x)$ à
la fonction nulle (au prix d'un Dirac à l'origine). Nous n'avons fait que
modifier notre distribution $D$ en lui ajoutant une mesure supportée sur
$\ZZ\setminus\{0\}$ donc elle reste une mesure supportée sur $\ZZ$ et sa
transformée de Fourier $E$ reste $1$-périodique et paire. Maintenant la
restriction de $E$ à $\OO{0,\frac12+\epsilon}$ est identiquement nulle. Donc $E$
restreinte à $\OO{-\frac12-\epsilon,\frac12+\epsilon}$ est une combinaison
linéaire finie de $\delta$, $\delta''$, \dots. Mais $E$ est $1$ périodique,
donc en fait
\[ E = P(\frac{d}{dx}) \sum_{k\in\ZZ} \delta_k(x)\] avec un polynôme $P$
pair. Ce qui signifie que $D$ qui est sa transformée de Fourier est de la
forme $R(x) \sum_{k\in\ZZ} \delta_k(x) = \sum_{k\in \ZZ} R(k)\delta_k(x)$ avec
$R$ un certain polynôme pair, qui doit en fait être nul à l'origine puisque
notre $D$ est nulle sur $\OO{-1,1}$.
Ainsi notre $f(s)$ du départ a
été réduit à une combinaison linéaire finie de $\zeta(s-2)$, $\zeta(s-4)$,
\dots

En conclusion, les hypothèses (1), (2), et (3) ont comme conséquence que la
série de Dirichlet de départ est une combinaison linéaire finie $f(s) =
\sum_{k\in\ZZ} c_k \zeta(s-2k)$ et le \autoref{th:1} est démontré.

\pdfbookmark[2]{Preuve du premier corollaire}{preuvecori}

Soit, pour $k\in \ZZ$, $f_k(s) = \zeta(s-2k)$ et $g_k(s) = \chi(s)
f_k(1-s)$. Par la relation $\Gamma(s+1)  = s\Gamma(s)$ et l'équation
fonctionnelle on a:
\begin{align*}
  g_k(s) &= \frac{s(s+1)\dotsm (s+2k-1)}{(-4\pi^2)^k} \zeta(s+2k)\qquad
  (k\ge0)\\
  g_k(s) &= \frac{(-4\pi^2)^{-k}}{(s+2k)(s+2k+1)\dotsm (s-1)} \zeta(s+2k)
  \qquad (k<0)
\end{align*}
Donc, pour tout $k$ donné, et tout $\sigma\gg0$ on a $g_k(s) \asymp
(-4\pi^2)^{-k} s^{2k}$ pour $\Re(s) = \sigma$, $|s|\to\infty$. Il en résulte
qu'une combinaison linéaire ne pourra être bornée, sur une droite quelconque
prise suffisamment à droite, que si elle ne comporte que des $k\le0$. Ceci
montre qu'en effet la condition (3${}'$) du Corollaire est équivalente à ce
que $f(s)$ soit une combinaison linéaire des $\zeta(s-2k)$, $k\le 0$. 

\pdfbookmark[2]{Preuve du second corollaire}{preuvecorii}

Et comme $g_k(\sigma)\sim (-4\pi^2)^{-k} \sigma^{2k}$ pour $\sigma\to+\infty$,
si on fait l'hypothèse \ref{cond:4} qu'une combinaison linéaire se rapproche
d'une limite $c$ plus vite que tout polynôme inverse ne se rapproche
de zéro, alors en particulier elle reste bornée lorsque $\sigma$ tend vers
$+\infty$ par valeurs réelles, ne peut comporter que des termes avec $k\le 0$,
et finalement que le seul terme avec $k=0$. Ainsi $f = c\zeta$.

\medskip

\pdfbookmark[2]{Preuve du \autoref{th:2}}{preuveii}

Nous débutons maintenant les preuves des Théorèmes \ref{th:2} et \ref{th:3}.
Soit $D_0(x) = \sum_{p\in\ZZ} \delta(x-p)$ la distribution de Poisson, qui est
sa propre transformée de Fourier (\emph{cf.}
\hyperref[annexe]{l'Annexe}). Considérons $\sum_{p\in\ZZ} e^{2\pi i (d+p)e}
{\delta(x- d - p)} = e^{2\pi i ex} D_0(x-d)$. Sa transformée de Fourier est
$e^{2\pi i de} e^{2\pi i dx}D_0(x+\nobreak e)$. Nous définirons donc
\[ D_{d,e}(x) = e^{-\pi i de}e^{2\pi i ex} D_0(x-d)\]
de sorte que $\cF$ agit comme la rotation d'angle $\frac\pi2$ dans le plan des
$(d,e)$:
\[ \cF(D_{d,e}) = D_{-e,d}\;,\]
Considérons  les distributions paires: $T_{d,e}(x) = \frac12
D_{d,e}(x) + \frac12 D_{d,e}(-x)$.  Comme $D_{d,e}(-x) = D_{-d,-e}(x)$, on a $T_{d,e}
= \frac12 D_{d,e} + \frac12 D_{-d,-e}$ et 
sa transformée de Fourier est:
\[ \cF(T_{d,e}) = T_{-e,d}\]
Compte tenu de
$T_{d+1,e} = e^{-\pi ie}T_{d,e}$, $T_{d,e+1} = e^{\pi i d} T_{d,e}$ et de $T_{-d,-e} =
T_{d,e}$, on peut toujours se ramener à $0<d\le 1$, $0<e\le1$, $0<d+e\le 1$,
et, pour $d+e=1$, par exemple $\frac12\le d \le 1$.

La transformée de Mellin $f(s)$ de la distribution paire $T_{d,e}$ est donnée
par: 
\[f(s) = \frac12 e^{-\pi i de} \sum_{x\equiv d\bmod 1,x>0} \frac{e^{2\pi i ex}}{x^{s}}
+  \frac12 e^{-\pi i de} \sum_{x\equiv -d\bmod 1,x>0} \frac{e^{-2\pi iex}}{x^{s}}\] La
fonction $g(s) = \chi(s) f(1-s)$ est:
\[ g(s) = \frac12 e^{\pi i de} \sum_{y\equiv e\bmod 1,y>0} \frac{e^{- 2\pi i
    dy}}{y^{s}} +  \frac12 e^{\pi i de} \sum_{y\equiv -e\bmod 1,y>0} \frac{e^{2\pi i d
    y}}{y^{s}}\] La fonction $\pi^{-\frac s2} \Gamma(\frac s2) f(s)$ a au plus
des pôles en $0$ et $1$: d'après \cite[Thm. 4.55]{jfcopoisson}, elle est
entière lorsque $d$ et $e$ ne sont pas entiers; lorsque $d\in \ZZ$ et
$e\notin\ZZ$, $\pi^{-\frac s2} \Gamma(\frac s2) f(s)$ possède un unique pôle,
qui est en $s=0$ et le résidu est $-e^{-\pi i de}$; lorsque $d\notin\ZZ$ et
$e\in\ZZ$, un unique pôle qui est en $s=1$ et de résidu $e^{\pi i de}$;
lorsqu'à la fois $d$ et $e$ sont entiers, un pôle en $s=0$ de résidu
$-(-1)^{de}$ et un pôle en $s=1$ de résidu $(-1)^{de}$.

Comme pour le \autoref{th:1} les hypothèses faites dans les Théorèmes
\ref{th:2} ou \ref{th:3} assurent que $f$ et $g$ sont des fonctions modérées,
et donc transformées de Mellin respectives de distributions paires $D$ et $E_1$ identiquement nulles dans un voisinage de l'origine, et telles que $E
= \cF(D)$ ne diffère de $E_1$ que par une distribution quasi-homogène (paire)
$Q$. Pour le \autoref{th:2}, on a donc:
\[ D = \sum_{n\ge1} a_{n} (\delta_{x_n} + \delta_{-x_n})\]
\[ E = Q + \sum_{n\ge1} b_{n} (\delta_{y_n} + \delta_{-y_n})\] Supposons que
l'on puisse trouver un nombre fini de suites arithmétiques de raison $1$
recouvrant le support de $D$ pour $x>0$. On a alors un nombre fini de nombres
réels $d_1$, \dots, $d_N$, distincts modulo $1$, obtenus comme les
représentants modulo $1$ dans $\OF{-\frac12,\frac12}$ de tous les $\pm
x_n$, $n\ge1$ ($a_n\neq0$), et qui sont tels que le support de $D$ dans $\RR$
est recouvert par les $d_j+\ZZ$. Les hypothèses ne font pas que les $d_j$
soient parmi les $\pm x_n$ mais assurent que chaque $d_j$ est tel que soit
$d_j+K$, soit $-d_j+K$ est un $x_n$ lorsque $K$ est suffisamment grand (et
réciproquement pour $n$ avec $a_n\neq0$).

Soit $\phi(x) = \prod_{1\le j\le N} (e^{2\pi i x} - e^{2\pi i
    d_j}) = \sum_{k=0}^N c_k e^{2\pi i (N-k)x}$, $c_0 = 1$, $c_N = \pm1$. Comme $D$ est une mesure on a:
  \begin{gather*}
     \phi(x) D(x) = 0\qquad\text{et donc}\qquad 
     E(x+N) + c_1 E(x+ N-1) + \dots + c_N E(x) = 0\;.
  \end{gather*}
Pour $x>0$ l'expression
\[ Q(x+N) + c_1 Q(x+N-1) + \dots + c_N Q(x)\] est une fonction analytique
(dans $\CC\setminus\OF{-\infty,0}$). Or elle a aussi, par l'équation pour $E$,
un support discret, et est par suite identiquement nulle. Et par conséquent
aussi \[ x>0\implies E_1(x+N) + c_1 E_1(x+ N-1) + \dots + c_N E_1(x) = 0\;.\]
Montrons aussi au passage que $Q$ est (pour $x>0$) constante. Comme $|c_N| =
1$, la formule de récurrence permet de voir que $Q$ est analytique dans
$\CC\setminus\OF{-\infty,-1}$, puis dans $\CC\setminus\OF{-\infty,-2}$,
\dots{} et donc que $Q$ est une fonction entière. Comme elle a a priori
croissance polynomiale à l'infini, c'est qu'elle est un polynôme. Le degré
d'un polynôme $Q(x)$ ne change pas par son remplacement par $Q(x+1)- \omega
Q(x)$, lorsque $\omega\neq1$, donc l'équation pour $Q$ donne $Q=0$ si $0$
n'est pas dans la liste des $d_j$, et se réduit à $Q(x+1) - Q(x) = 0$ donc à
$Q$ constante sinon. Ceci est valable pour $x>0$ et par parité la distribution
$Q$ est une constante sur $\RR\setminus\{0\}$.

En ce qui concerne $E_1$, la formule montre que si $x>N$ appartient à son
support alors il en est de même de l'un de $x-1$, $x-2$, \dots, $x-N$. Comme
le support de $E_1$ est discret, cela prouve (comme affirmé dans l'énoncé)
qu'en effet il est contenu dans un nombre
fini de suites arithmétiques de raison $1$. 

La situation est maintenant symétrique entre $D$ et $E$ et compte tenu de ce
que nous avons déjà établi pour la restriction de $Q$ à $x\neq 0$, nous en
déduisons que $\cF(Q)$ est constant pour $x\neq0$ et donc au final que $Q$
est de la forme $-\alpha + \beta \delta$. Quitte à remplacer $D$ par $D+\alpha
\delta$ et $E$ par $E +\alpha$ on peut dorénavant supposer que $D$ comme $E$
ont, au plus, chacun un Dirac à l'origine. Cette modification éventuelle de
$D$ (qui n'en est pas une de la série de Dirichlet d'origine) peut nous
amener à joindre $0$ à la liste des $d_j$ (et à remplacer $N$ par
$N+1$). Dorénavant $E$ est une mesure paire, vérifiant une récurrence
\[ \prod_{j} (\tau - \exp(2\pi i d_j)) E = 0\qquad (\tau(E)(x) = E(x+1))\] et
supportée sur un nombre fini de translatés de $\ZZ$. Soient $e_1$, \dots,
$e_M$ les représentants dans $]-\frac12,+\frac12]$ de ce support de $E$ pris
modulo $1$. On a:
\[ E(x) = \sum_{1\le k\le M}\, \sum_{p\in\ZZ} c_k(p) \delta(x - e_k - p) \]
Chacune des $M$ suites $(c_k(p))_{p\in\ZZ}$ doit vérifier la récurrence
linéaire
\[ \prod_{1\le j\le N} (\tau - \exp(2\pi i d_j)) c_k = 0\qquad (\tau(c_k)(p) =
c_k(p+1))\] dont l'espace vectoriel des solutions est de dimension $N$ et
engendré par les fonctions linéairement indépendantes $p\mapsto \exp(2\pi i
d_j p)$, $1\le j \le N$. Cela prouve que $E$ est une combinaison linéaire des
$NM$ distributions $e^{2\pi i d_j x} D_0(x- e_k)$ considérées précédemment
(avec $D_0$ la distribution de Poisson), $1\le j \le N$, $1\le k\le M$. Mais
comme $E$ est une distribution paire elle est aussi combinaison linéaire des
parties paires de ces distributions, donc des $T_{e_k,d_j}$. Par conséquent la
distribution originale $D$ (incorporant, comme nous l'avons indiqué, un
certain multiple du Dirac à l'origine) est une combinaison linéaire des
$T_{-d_j,e_k}$, ou encore des $T_{d_j,e_k}$ puisque l'ensemble des $d_j$ est
stable par passage à l'opposé (sauf pour un éventuel $\frac12$, mais $\frac12
= -\frac12+1$ et $T_{d+1,e} = e^{-\pi i e}T_{d,e}$).
Le \autoref{th:2} est  prouvé.

\pdfbookmark[2]{Preuve du \autoref{th:3}}{preuveiii}

Je passe maintenant à la preuve du \autoref{th:3}. En exploitant les
hypothèses comme précédemment nous avons ici une distribution tempérée paire
$D = \sum_{n\ge1} a_n(\delta_{x_n} + \delta_{-x_n})$, dont le support est
inclus (par hypothèse) dans un nombre fini $N$ de translatés de $\ZZ$, et la
transformée de Fourier $E$ est de la forme $Q + E_1$ avec $Q$ quasi-homogène,
$E_1$ paire, égale dans un certain intervalle ouvert $\OO{-Y_0,Y_0}$ ($Y_0>N$)
à $\sum_{1\le n \le M} b_n(\delta_{y_n} + \delta_{-y_n})$. Comme précédemment,
par ce que nous savons du support de la mesure $D$ nous obtenons une certaine
récurrence linéaire $\prod_{1\le j \le N} (\tau - \exp(2\pi i d_j)) E =
0$. Prenons $\epsilon>0$ avec $N+\epsilon < Y_0$ et de sorte que l'intervalle
$]0,\epsilon[$ ne rencontre modulo $1$ aucun des $\pm y_n$, $1\le n \le
M$. Sur cet intervalle $E_1(x)$, \dots, $E_1(x+N)$ sont identiquement nuls, et
ainsi $\prod_{1\le j \le N} (\tau - \exp(2\pi i d_j)) Q = 0$. On peut donc
exprimer $Q(x)$ en fonction de $Q(x+1)$, \dots, $Q(x+N)$, et par analyticité
cela restera valable pour $x = z \in \CC\setminus\OF{-\infty,0}$, puis permet
d'étendre le domaine d'analyticité de $Q(z)$ à $\CC\setminus\OF{-\infty,-1}$,
etc\dots{}, donc $Q(z)$ est une fonction entière et finalement un
polynôme. Pour les $\omega\neq1$ la transformation $P(x)\mapsto P(x+1)- \omega
P(x)$ conserve le degré des polynômes, donc notre polynôme $Q(x)$ ou est
directement vu comme étant nul, ou vérifie $Q(x+1) = Q(x)$ et est
constant. Par parité nous avons ainsi que la restriction de $Q$ à $x\neq0$ est
une constante. Nous pouvons alors soustraire un Dirac à l'origine à $D$ afin
de retirer cette constante à $E$, en maintenant la relation $\cF(D) = E$. La
restriction de $E$ à $\OO{-\epsilon,+\epsilon}$ est paire et supportée à
l'origine donc une combinaison de $\delta$, $\delta''$, \dots{}. Cependant $E$
n'a en $1$, $2$, \dots, $N$ au plus que des Dirac, et c'est ainsi aussi le cas
à l'origine, par la récurrence. La restriction de $E$ à
$\OO{-\epsilon,N+\epsilon}$ est ainsi une somme d'un nombre fini de Dirac et
par la relation de récurrence, il en résulte que $E$ est elle-même une mesure,
au support discret inclus dans un nombre fini de suites arithmétiques de raison
$1$. À ce stade nous nous sommes ramenés aux hypothèses du Théorème précédent.

\section{Annexe:  preuve de l'identité de Poisson distributionnelle}
\label{annexe}

Soit $D(x) = \sum_{n\in \ZZ} \delta(x-n)$ la distribution de Poisson, et $E$
sa transformée de Fourier, au sens des distributions tempérées. De $D(x+1) =
D(x)$ il résulte $(e^{2\pi i x} - 1) E(x) = 0$ ce qui établit que $E$ est une
mesure\footnote{puisque les seules distributions annulées par $x$ sont le
  Dirac et ses multiples!} supportée sur $\ZZ$. De $(e^{2\pi i x} - 1) D(x) =
0$ il résulte que $E$ est $1$-périodique. Donc $E = cD$ pour une certaine
constante $c$.

Soit $f$ une fonction paire, infiniment dérivable à support compact, positive,
non identiquement nulle. Soit $k = f\mathbin{\mbox{\raisebox{-.25\height}{*}}}f + \cF(f)^2$. La fonction de Schwartz $k$
est paire, positive, et sa propre transformée de Fourier. Et $k(0)>0$. Donc
$0 < (D,k) = (E,k)$ et $c=1$.

\end{document}